\documentclass[12pt]{article}
\usepackage{amsfonts}
\usepackage{amsthm,amsmath, amssymb}
\usepackage{graphicx}
\usepackage[authoryear]{natbib}
\usepackage{verbatim}
\usepackage{bbm}
\usepackage{authblk}
\usepackage{setspace}
\newtheorem{theorem}{Theorem}
\newtheorem{definition}{Definition}
\newtheorem{corollary}{Corollary}
\newtheorem{lemma}{Lemma}
\newtheorem{remark}{Remark}

\providecommand{\bo}[1]{\boldsymbol{#1}}

\newcommand{\bmat}{\begin{pmatrix}}
\newcommand{\emat}{\end{pmatrix}}

\usepackage{setspace}
\title{Separation of uncorrelated stationary time series using autocovariance matrices}
\author[1]{Jari Miettinen\thanks{\it Address for correspondence: Jari Miettinen, Department of Mathematics and Statistics, 40014 University of Jyv\"askyl\"a, Finland. E-mail: jari.p.miettinen@jyu.fi}}
\author[2]{Katrin Illner}
\author[3]{Klaus Nordhausen}
\author[3]{Hannu Oja}
\author[1]{Sara Taskinen}
\author[4]{Fabian J. Theis}
\affil[1]{\small University of Jyv\"askyl\"a}
\affil[2]{\small Helmholtz Center Munich, Germany}
\affil[3]{\small University of Turku, Finland.}
\affil[4]{\small Technical University of Munich, Germany}

\begin{document}

\maketitle

\begin{abstract}
Blind source separation (BSS) is a signal processing tool, which is widely used in various fields.
Examples include biomedical signal separation, brain imaging and economic time series applications.
In BSS, one assumes that the observed $p$ time series are linear combinations of $p$ latent uncorrelated
weakly stationary time series. The aim is then to find an estimate for an unmixing matrix, which transforms
the observed time series back to uncorrelated latent time series. In SOBI (Second Order Blind Identification)
joint diagonalization of the covariance matrix and autocovariance matrices with several lags is used to
estimate the unmixing matrix. The rows of an unmixing matrix can be derived either one by one (deflation-based
approach) or simultaneously (symmetric approach). The latter of these approaches is well-known especially in signal
processing literature, however, the rigorous analysis of its statistical properties has been missing so far.
In this paper, we fill this gap and investigate the statistical properties of the symmetric SOBI estimate
in detail and find its limiting distribution under general conditions. The asymptotical efficiencies of
symmetric SOBI estimate are compared to those of recently introduced deflation-based SOBI estimate under
general multivariate MA$(\infty)$ processes. The theory is illustrated by some finite-sample simulation
studies as well as a real EEG data example.

Keywords: Asymptotic normality; Blind source separation; Joint diagonalization; MA($\infty$); SOBI
\end{abstract}

%\begin{keyword}[class=AMS]
%\kwd[Primary ]{62H05, 62H10} \kwd[; secondary]{62H12}
%\end{keyword}

%\keywords{Asymptotic normality; Blind source separation; Joint diagonalization; MA($\infty$); SOBI}

\section{Introduction}
In blind signal separation or blind source separation (BSS) one assumes that a $p$-variate
observable random vector $\bo x$ is a linear mixture
of $p$-variate latent source vector $\bo z$. The model can thus be written as
%\begin{align}
%\label{BSSmodel0}
$\bo x=\bo\mu+\bo\Omega\bo z$,
%\end{align}
where the full rank $p \times p$ matrix $\bo\Omega$ is so called {\it mixing matrix} and $\bo z$ is
a random $p$-vector with certain preassigned properties. The $p$-vector $\bo\mu$ is a
location parameter and usually considered as a nuisance parameter, since the main goal in BSS
is to find an estimate for an {\it  unmixing matrix} $\bo\Gamma$ such that $ \bo\Gamma \bo x\sim \bo z$,
based on a $p\times n$ data matrix $\bo X=(\bo x_1,\dots,\bo x_n)$ from the distribution of $\bo x$.

The BSS model was formulated for signal processing  and computer science applications in the early
1980's and since then, many approaches have been suggested to solve the problem under various assumptions
on $\bo z$. For an account of the early history of BSS, see \cite{JuttenTaleb:2000}. The most popular
BSS approach is {\it independent component analysis (ICA)}, which assumes that $E(\bo z)=0$ and
$E(\bo z\bo z')=\bo I_p$, and that the components of $\bo z$ are mutually independent.
The ICA model is a semiparametric model as the marginal distributions of the components of $\bo z$ are
not specified at all. For identifiability of the parameters, one has to assume, however, that at most
one of the components is normally distributed.
Typical algorithms for ICA use centering and  whitening as preprocessing steps: Write 
$\bo\Omega=\bo U\bo\Lambda \bo V'$ for the singular value decomposition (SVD) of the mixing matrix
$\bo\Omega$. Then, under the above mentioned assumptions on $\bo z$, $E(\bo x)=\bo\mu$ and
$Cov(\bo x)=\bo\Sigma=\bo U\bo\Lambda^2 \bo U'$, and therefore $\bo V\bo U'\bo\Sigma^{-1/2}(\bo x-\bo\mu)=\bo z$.
An iterative algorithm can then be applied to $\bo\Sigma^{-1/2}(\bo x-\bo\mu)$ to find the orthogonal 
matrix $\bo V\bo U'$. For an overview of ICA from a signal processing perspective, 
see for example~\cite{HyvarinenKarhunenOja:2001} and~\cite{ComonJutten:2010}.

Since the late 1990's, there has been an increasing interest in ICA methods among statisticians.
\citet{OjaSirkiaEriksson:2006}, for example, used two scatter matrices, $\bo S_1$ and $\bo S_2$, with the
independence property to solve the ICA problem. We say that a $p\times p$ matrix valued functional $\bo S(F)$ is a
{\it scatter matrix} if it is symmetric, positive definite and affine equivariant in the sense that
$\bo S(F_{\bo A\bo x+\bo b})=\bo A\bo S(F_{\bo x})\bo A'$ for all full-rank $p\times p$ matrices $\bo A$ and
for all $p$-vectors $\bo b$. Moreover, a scatter matrix $\bo S(F)$ has the {\it independence property}
if $\bo S(F_{\bo z})$ is a diagonal matrix for all $\bo z$ with independent components.
An unmixing matrix $\bo\Gamma$ and a diagonal matrix $\bo\Lambda$ (diagonal elements in a decreasing order)
 then solve the estimating equations
\[
\bo\Gamma \bo S_1\bo\Gamma'=\bo I_p\ \ \mbox{and}\ \ \bo\Gamma \bo S_2\bo\Gamma'=\bo \Lambda.
\]
The independent components in $\bo\Gamma\bo x$ are thus standardized with respect to $\bo S_1$ and uncorrelated with
respect to $\bo S_2$, and $\bo\Gamma$ and $\bo\Lambda$ can be found as eigenvector-eigenvalue solutions
for  $\bo S_1^{-1}\bo S_2$. An example of such ICA methods is the classical FOBI (fourth order blind identification)
method which uses $\bo S_1(F_{\bo x})=Cov(\bo x)$ and
$$
\bo S_2(F_{\bo x})= E\left[(\bo x-E(\bo x)) (\bo x-E(\bo x))'Cov(\bo x)^{-1} (\bo x-E(\bo x))(\bo x-E(\bo x))' \right]  .
$$
See~\citet{OjaSirkiaEriksson:2006}, for example. 
The unmixing matrix estimate $\hat{\bo\Gamma}$ is naturally obtained by replacing the population 
values by their sample counterparts. 

The limiting statistical properties of several ICA unmixing matrix estimates have been developed
quite recently and mainly for iid data: For the limiting behavior of the unmixing matrix estimate based on 
two scatter matrices, see \cite{IlmonenNevalainenOja:2010a}. For other recent work and different estimation 
procedures for ICA, see for example \cite{HastieTibshirani:2003}, \cite{ChenBickel:2006}, \cite{BonhommeRobin:2009},
\cite{IlmonenPaindaveine:2011}, \cite{AllassonniereYounes:2012}, \cite{SamworthYuan:2013} and
\cite{HallinMehta:2013}. 

In applications such as the analysis of medical
images or signals (EEG, MEG or fMRI) and  financial or geostatistical times series, the assumption
of independent observations does not usually hold. Nevertheless, ICA has been considered in this context 
in~\cite{ChenHardleSpokoiny:2007},~\cite{GarciaFerrerGonzalezPrietoPena:2011},~\cite{GarciaFerrerGonzalezPrietoPena:2012},~\cite{LeeShenTruongLewisHuang:2011},~\cite{Poncela:2012} and~\cite{Schachtneretal:2008}
among others. See also \cite{MattesonTsay:2011} for a slightly different model. Apart from these results, 
other BSS models have been developed for time series  data in signal processing literature. 
The so called  AMUSE (Algorithm for Multiple Unknown Signals Extraction) and SOBI (Second Order
Blind Identification) procedures for the stationary time series BSS models were suggested 
by~\cite{TongSoonHuangLiu:1990} and ~\cite{BelouchraniAbedMeraimCardosoMoulines:1997}, respectively.
\cite{MiettinenNordhausenOjaTaskinen:2012,MiettinenNordhausenOjaTaskinen:2014} provided careful
analysis of the statistical properties of the AMUSE and so called deflation-based SOBI  estimates. 
\cite{Nordhausen:2014} considered methods that assume only local stationarity.

In this paper we will continue the work of 
\cite{MiettinenNordhausenOjaTaskinen:2012,MiettinenNordhausenOjaTaskinen:2014} and derive the 
statistical properties of so called symmetric SOBI estimate. We will use a real EEG data example to illustrate how
the theoretical results derived in this paper may be used to measure the accuracy of the unmixing matrix 
estimates. The structure of this paper is as follows. 
In Section~\ref{SOSmodel} we introduce the blind source separation model which assumes second order 
stationary components. In Section~\ref{SOSfunctionals} we recall the definition for the deflation-based SOBI 
functional and  define symmetric SOBI functional using the Lagrange  multiplier technique. 
The theoretical properties of deflation-based and symmetric SOBI estimators 
are given in general case in Section~\ref{properties}. Further, in Section~\ref{MAinf} the limiting 
distributions of the two SOBI estimators will be more concretely compared under the assumption of 
$MA(\infty)$ processes. The theoretical results are illustrated using simulation studies in 
Section~\ref{mdsec} and a real EEG data example in Section~\ref{EEG} before the paper is concluded in 
Section~\ref{conclusion}. Asymptotical results for the symmetric SOBI estimates are proven in the appendix.

\section{Second order source separation model}\label{SOSmodel}

We  assume that the observable $p$-variate time series $\bo x=(\bo x_t)_{t=0,\pm 1,\pm 2, \ldots}$ are distributed
according to 
\begin{align}
\label{BSSmodel}
\bo x_t=\bo\mu+\bo\Omega\bo z_t,\ \ t=0,\pm 1,\pm 2,\ldots,
\end{align}
where  $\bo\mu$ is a $p$-vector, $\bo\Omega$ is a full-rank $p\times p$ { mixing matrix} and
$\bo z=(\bo z_t)_{t=0,\pm 1,\pm 2, \ldots}$ is a $p$-variate latent  time series that  satisfies
\begin{itemize}
\item[(A1)] $E(\bo z_t)=0$ and  $E(\bo z_t\bo z_t')=\bo I_p$.
\item[(A2)] $E(\bo z_t\bo z_{t+\tau}')= E(\bo z_{t+\tau}\bo z_t')= \bo\Lambda_\tau$ is diagonal for all $\tau=1,2,\ldots$.
\end{itemize}

\begin{figure}[htb]
\centering
\includegraphics[width=0.8\textwidth]{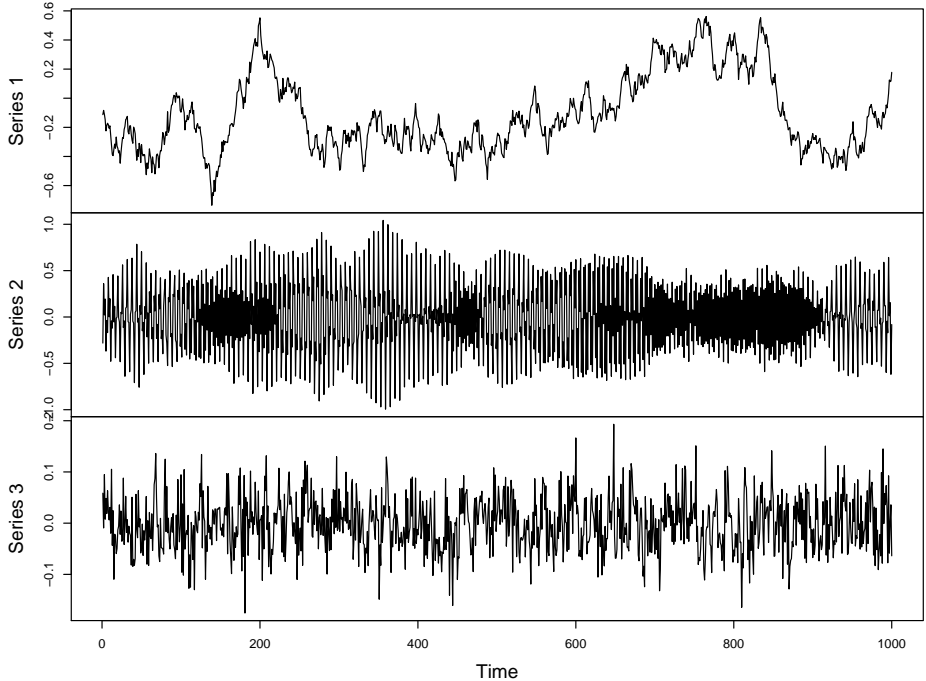}
%\hskip-0.5truecm
%\includegraphics[width=7truecm,height=7truecm]{.pdf}
\caption{\label{sources} Example time series $\bo z$: Three independent stationary AR series. }
\end{figure}

This is again a semiparametric model as only the moment assumptions (A1)-(A2) of the time series in $\bo z$ are made.
The assumptions state that the $p$ time series in $\bo z$ are weakly stationary and uncorrelated. This model is
called the {\it second order source separation (SOS) model}. A model with stronger assumptions, i.e. the
{\it independent component time series model}, is obtained if the condition (A2) is replaced by the condition
\begin{itemize}
\item[(A2*)] the $p$ times series in $\bo z$ are mutually independent and
$E(\bo z_t\bo z_{t+\tau}')= E(\bo z_{t+\tau}\bo z_t')= \bo\Lambda_\tau$ is diagonal for all $\tau=1,2,\ldots$.
\end{itemize}

\begin{figure}[htb]
\centering
\includegraphics[width=0.8\textwidth]{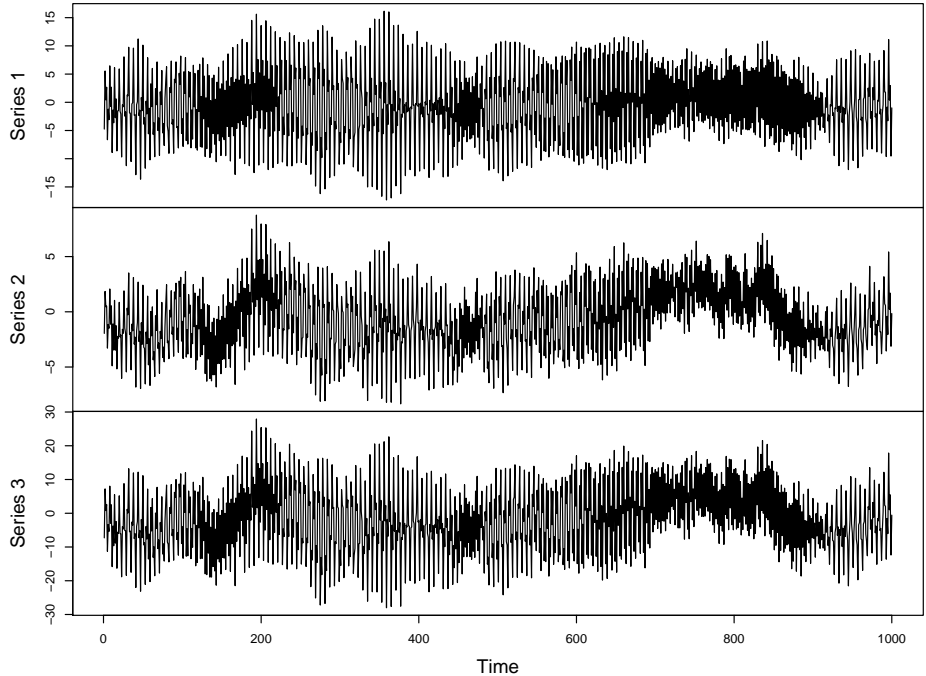}
\caption{\label{data}Time series $\bo x$: Mixtures of the three independent stationary AR series in
Figure~\ref{sources}. }
\end{figure}

Figure~\ref{sources} serves as an example of a 3-variate time series $\bo z$ with three independent
components, namely  $AR$ processes with coefficient vectors $(0.9, 0.09)$, $(0,0,-0.99)$ and $(0,0.3)$.
The observable 3-variate time series $\bo x$ consisting of three different mixtures of the latent time
series in $\bo z$ are shown in Figure~\ref{data}. Given the observed time series
$(\bo x_1,\ldots,\bo x_T)$, the aim is to find an estimate $\hat{\bo\Gamma}$ of an unmixing matrix
$\bo\Gamma$ such that $\bo\Gamma\bo x$ has uncorrelated components. Clearly, $\bo\Gamma=\bo C\bo\Omega^{-1}$ is an
unmixing matrix for any $p\times p$ matrix $\bo C$ with exactly one nonzero element in each row and in each
column. Notice that the signs and order of the components of $\bo z$  and the signs and order of the
columns  of $\bo\Omega$  are confounded also in the BSS model. Additional assumptions are therefore needed
in order to study the consistency and asymptotical properties of $\hat{\bo\Gamma}$. 
Contrary to ICA in the iid case, the mixing matrix may now be identifiable for any number of 
gaussian components. However, as we will see later in this paper, weak assumptions on the autocovariance 
matrices $\bo\Lambda_\tau$, $\tau=1,2,\dots$, have to be made for the identifiability of our functionals and 
for the study of the asymptotic properties of corresponding estimates.

\section{BSS functionals based on autocovariance matrices} \label{SOSfunctionals}

\subsection{Joint diagonalization of autocovariance matrices}
\label{amuse}

\begin{figure}[htb]
\centering
\includegraphics[width=0.8\textwidth]{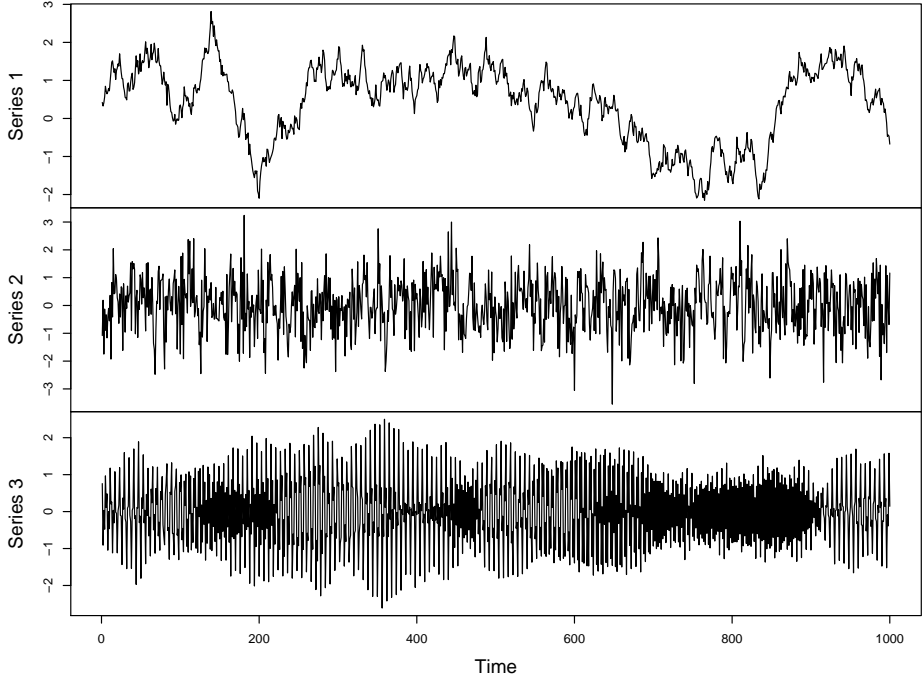}
\caption{Time series $\hat{\bo\Gamma}_1\bo x$ obtained with AMUSE from $\bo x$ in
Figure~\ref{data}. }
\label{AMUSEfig}
\end{figure}

The separation of uncorrelated stationary time series can be solely based on autocovariances and
cross-autocovariances of the $p$ time series. Assume that $\bo x$ follows a centered SOS model with
$\bo\mu=\bo 0$. This is not a restriction in our case as
the asymptotical properties of the estimated autocovariances are the same for known and unknown $\bo\mu$.
It then follows that
$$
 E(\bo x_t\bo x'_{t+\tau})=\bo\Omega \bo\Lambda_\tau \bo\Omega',\ \ \ \tau=0,1,2,\dots.
$$
Assume also that,
for some lag $\tau>0$, the diagonal elements of the autocovariance matrix
$E(\bo z_t\bo z'_{t+\tau})=\bo\Lambda_\tau$ are distinct.  An unmixing matrix functional $\bo\Gamma_\tau$  is
then defined as a $p\times p$ matrix that satisfies
$$
\bo\Gamma_\tau   E(\bo x_t\bo x'_{t})  \bo\Gamma_\tau'=\bo I_p\ \ \ \text{and}\ \ \ \bo\Gamma_\tau  E(\bo x_t\bo x'_{t+\tau})
\bo\Gamma_\tau'=\bo P_\tau \bo\Lambda_\tau \bo P_\tau',
$$
where  $\bo P_\tau \bo\Lambda_\tau \bo P_\tau' $ is a diagonal matrix with the same diagonal
elements as in $\bo\Lambda_\tau$ but in a decreasing order. (As in PCA and ICA, the signs of the rows
of $\bo\Gamma_\tau$ are not fixed in this definition.)
The components of  $\bo\Gamma_\tau \bo x$ are thus the components of $\bo z$
in a permuted order.  The permutation matrix $\bo P_\tau$ remains unidentifiable. Notice that $\bo\Gamma_\tau$
is affine equivariant, that is, the transformation $\bo x \to \bo A\bo x$ with a full-rank
$p\times p$ matrix $\bo A$ induces the transformation $\bo\Gamma_\tau\to \bo\Gamma_\tau \bo A^{-1}$.
This implies that $\bo\Gamma_\tau\bo x$ does not depend on the mixing matrix $\bo\Omega$ at all.
The corresponding sample statistic, the so called  AMUSE (Algorithm for Multiple Unknown Signals
Extraction) estimator, was proposed by~\cite{TongSoonHuangLiu:1990}. Figure~\ref{AMUSEfig} shows the
estimated latent sources obtained with AMUSE, $\tau=1$, from the data in Figure~\ref{data}.
See \cite{MiettinenNordhausenOjaTaskinen:2012} for a recent study of the statistical properties of
the AMUSE estimate.

The drawback of the AMUSE procedure is the assumption that, for the chosen lag $\tau$, the
eigenvalues in $\bo\Lambda_\tau$ must be distinct. This is of course never known in practice.  
Therefore, the choice of $\tau$ may have a huge impact on the performance of the method, as 
only information coming from $S_0=E(\bo x_t\bo x'_{t})$ and $S_\tau=E(\bo x_t\bo x'_{t+\tau})$ is used.
To overcome this drawback,~\cite{BelouchraniAbedMeraimCardosoMoulines:1997} proposed the SOBI
(Second Order Blind Indentification) algorithm that aims to jointly diagonalize several
autocovariance matrices as follows. Let $\bo S_1,\dots,\bo S_K$ be $K$ autocovariance matrices with
distinct lags $\tau_1,\dots,\tau_K$. The $p\times p$ unmixing matrix functional
$\bo\Gamma=(\bo{\gamma}_1,\dots,\bo{\gamma}_p)'$ is the matrix that minimizes
$$
\sum_{k=1}^K||\mbox{off}(\bo\Gamma \bo S_{k}\bo\Gamma')||^2
$$
under the constraint $\bo\Gamma \bo S_0\bo\Gamma'=\bo I_p$, or, equivalently, maximizes
\begin{align*}
%\label{max1}
\sum_{k=1}^K||\mbox{diag}(\bo\Gamma \bo S_{k}\bo\Gamma')||^2
=\sum_{j=1}^p\sum_{k=1}^K(\bo\gamma_j' \bo S_{k}\bo\gamma_j)^2
\end{align*}
under the same constraint. Here we write diag$(\bo S)$ for a $p\times p$ diagonal matrix with the diagonal
elements as in $\bo S$ and off$(\bo S)=\bo S-$diag$(\bo S)$.

Next notice that, as $\bo\Gamma \bo S_0\bo\Gamma'=\bo I_p$, then $\bo\Gamma=\bo U\bo S_0^{-1/2}$ for some orthogonal
$p\times p$ matrix $\bo U=(\bo u_1,\dots,\bo u_p)'$. If then
$\bo R_k=\bo S_0^{-1/2} \bo S_k \bo S_0^{-1/2}$, $k=1,\dots,K$,
are the autocorrelation matrices, the solution for $\bo U$ can be found by maximizing
\begin{align}
\label{max1b}
\sum_{k=1}^K||\mbox{diag}(\bo U \bo R_{k} \bo U')||^2
=\sum_{j=1}^p\sum_{k=1}^K(\bo u_j' \bo R_{k}\bo u_j)^2
\end{align}
under the orthogonality constraints $\bo U\bo U'=\bo I_p$.

In the literature, several algorithms to solve the maximization problem in~(\ref{max1b}) are proposed:
In deflation-based approach, the rows of $\bo U$ are found one by one using some pre-assigned rule.
In the symmetric approach, the rows are found simultaneously. The solution $\bo\Gamma=\bo U\bo S_0^{-1/2}$
naturally depends on the approach as well as on the concrete algorithm used in the optimization. In the 
following we consider  deflation-based and symmetric approaches in more detail. 
%To treat the maximization problem, we use in the following the functions
%\[
%\label{functions}
%D(\bo u)=\sum_{k=1}^K (\bo u'R_k \bo u)^2
%\ \ \mbox{and}\ \
%T(\bo u)=\frac 12 \frac {\partial}{\partial \bo u} D(\bo u)=\sum_{k=1}^K (\bo u'R_k \bo u) R_k\bo u.
%\]

\subsection{Deflation-based approach}
\label{sobid}

In the deflation-based approach, the rows of an unmixing matrix functional
$\bo\Gamma=(\bo\gamma_1,\dots,\bo\gamma_p)'$ are found one by one so that $\bo\gamma_j$, $j=1,\dots,p-1$,
maximizes
\begin{align}
\label{max2}
\sum_{k=1}^K(\bo\gamma_j'\bo S_k\bo\gamma_j)^2,
\end{align}
under the constraints $\bo\gamma'_i\bo S_0\bo\gamma_j=\delta_{ij}$, $i=1,\dots,j$.
Recall that the Kronecker delta $\delta_{ij}=1\ (0)$ as $i=j\ (i\ne j)$.

The solution  $\bo\gamma_j$ then optimizes the Lagrangian function
$$
L(\bo\gamma_j,\bo\theta_j)=\sum_{k=1}^K (\bo\gamma_j'\bo S_k\bo\gamma_j)^2
-\theta_{jj}(\bo\gamma'_j\bo S_0\bo\gamma_j-1)-\sum_{i=1}^{j-1}\theta_{ji}\bo\gamma'_i\bo S_0\bo\gamma_j,
$$
where $\bo\theta_j=(\theta_{j1},\dots,\theta_{jj})'$ are the Lagrangian multipliers.
Write
\begin{align}\label{T}
\bo T(\bo{\gamma})=\sum_{k=1}^K(\bo{\gamma}'\bo S_k\bo{\gamma})\bo S_k\bo{\gamma}.
\end{align}
The unmixing matrix functional $\bo\Gamma$ found in this way then satisfies the following estimating
equations~\citep{MiettinenNordhausenOjaTaskinen:2014}.

\begin{definition}
\label{djd}
The deflation-based unmixing matrix functional
$\bo\Gamma=(\bo{\gamma}_1,\dots,\bo{\gamma}_p)'$
solves the $p-1$ estimating equations
$$
\bo{T}(\bo{\gamma}_j)=\bo S_0 \left(\sum_{r=1}^j\bo{\gamma}_r\bo{\gamma}'_r \right)\bo{T}(\bo{\gamma}_j),\
\ \ j=1,\dots,p-1.
$$
\end{definition}

Recall that $\bo\Gamma=\bo U\bo S_0^{-1/2}$ with some orthogonal matrix $\bo U=(\bo u_1,\dots,\bo u_p)'$
and, in the deflation-based approach,  the rows of $\bo U$ are found one by one as well. The estimating
equations then suggest the following fixed point algorithm for the deflation-based solution.
After finding ${\bo u}_1,\dots,{\bo u}_{j-1}$, the following two steps are repeated until convergence
to get $\bo u_j$.
\begin{align*}
&\text{step 1: } \ \   {\bo u}_{j} \leftarrow \left(\bo I_p-\sum_{i=1}^{j-1}{\bo u}_i{\bo u}'_i\right) {\bo T}({\bo u}_{j}). \\
&\text{step 2:}\ \ {\bo u}_{j} \leftarrow||{\bo u}_{j}||^{-1}  {\bo u}_{j}.
\end{align*}
Here $\bo T(\bo u)=\sum_{k=1}^K (\bo u'\bo R_k \bo u) \bo R_k\bo u$. Notice that the algorithm naturally needs
initial values for each ${\bo u}_j$, $j=1,\dots,p-1$, and different initial values may change the
rows of the estimate and produce them in a permuted order. Therefore, for $\bo u_j$, one should use
several randomly selected initial values to guarantee that the true maximum in  (\ref{max2}) is
attained at each stage. For a more detailed study of this algorithm, see Appendix~A
in~\cite{MiettinenNordhausenOjaTaskinen:2014}.

\subsection{Symmetric approach}
\label{sobis}

In the symmetric approach, the rows of an unmixing matrix functional $\bo\Gamma=(\bo\gamma_1,\dots,\bo\gamma_p)'$
are found simultaneously. We then consider the maximization of
\begin{align*}
%\label{max3}
\sum_{j=1}^p\sum_{k=1}^K(\bo\gamma_j'\bo S_k\bo\gamma_j)^2,
\end{align*}
under the constraint $\bo\Gamma \bo S_0\bo\Gamma^T=\bo I_p$. The matrix $\bo\Gamma$ now optimizes the
Lagrangian function
$$
L(\bo\Gamma,\Theta)=\sum_{j=1}^p\sum_{k=1}^K (\bo\gamma_j'\bo S_k\bo\gamma_j)^2
-\sum_{j=1}^p \theta_{jj}(\bo\gamma'_j\bo S_0\bo\gamma_j-1)-
\sum_{j=1}^p\sum_{i=1}^{j-1} \theta_{ij}\bo\gamma'_i\bo S_0\bo\gamma_j,
$$
where the symmetric matrix $\bo\Theta=(\theta_{ij})$ contains the $p(p+1)/2$ Lagrangian multipliers of the
optimization problem.
%Write again
%$$
%\bo T(\bo \gamma)=\sum_{k=1}^K(\bo \gamma'\bo S_k\bo \gamma)\bo S_k\bo \gamma.
%$$
At the solution $\bo\Gamma$ we then have
\begin{align*}
2 \bo T(\bo \gamma_j)= \bo S_0 \left(2\theta_{jj} \bo \gamma_j+\sum_{i=1}^{j-1} \theta_{ij} \bo \gamma_i
+\sum_{i=j+1}^{p} \theta_{ji} \bo \gamma_i\right),
\end{align*}
where $\bo T(\bo\gamma)$ is as in~(\ref{T}). Multiplying both sides from the left by $\bo\gamma_i'$ gives
$
2\bo\gamma_i' \bo T(\bo\gamma_j)=\theta_{ij},
$
for  $i<j$,
and
$
2\bo\gamma_i' \bo T(\bo\gamma_j)=\theta_{ji},
$
for $i>j$.
Hence the solution $\bo\Gamma$ must satisfy the following estimating equations.
\begin{definition}
\label{sjd}
The symmetric unmixing matrix functional
$\bo\Gamma=(\bo{\gamma}_1,\dots,\bo{\gamma}_p)'$
solves the  estimating equations
\begin{equation*}
\bo\gamma_i'\bo T(\bo\gamma_j)=\bo\gamma_j'\bo T(\bo\gamma_i) \ \ \ \text{and}\ \ \
\bo\gamma_i'\bo S_0\bo\gamma_j=\delta_{ij},\ \ i,j=1,\dots,p.
\end{equation*}

\end{definition}

\medskip
\begin{figure}[htb]
\centering
\includegraphics[width=0.8\textwidth]{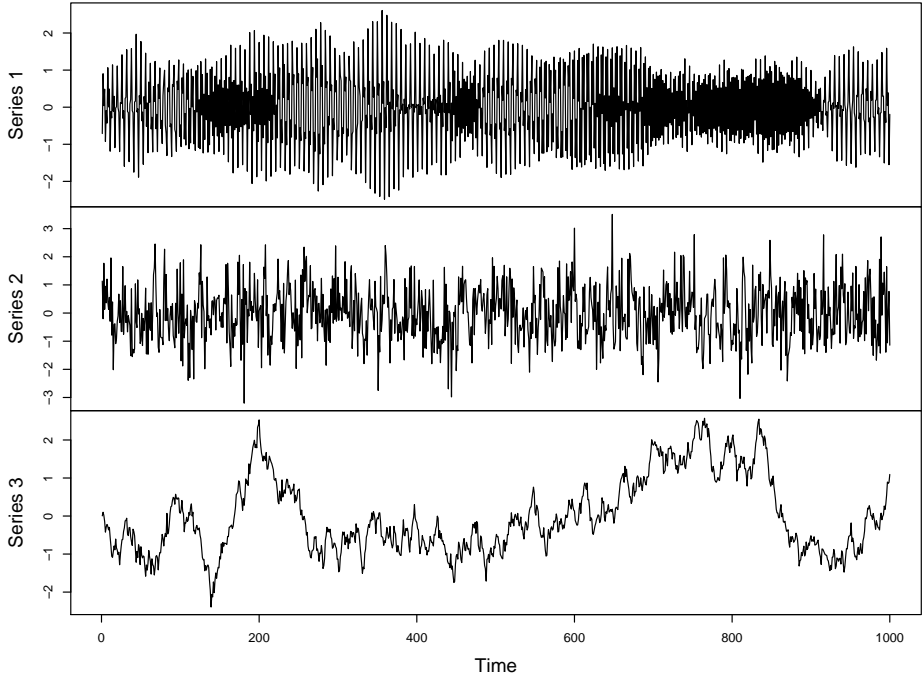}
\caption{\label{SOBIfig}Time series $\hat{\bo\Gamma}\bo x$ obtained with symmetric SOBI from $\bo x$ in
Figure~\ref{data}. }

\end{figure}

Notice that the exact joint diagonalization is possible only if the matrices $\bo R_1,\dots,\bo R_K$ have
the same sets of eigenvectors. This is naturally true for the population matrices in the SOS model.
For estimated autocorrelation matrices from a continuous SOS model, the eigenvectors are however almost
surely different. As $\bo\Gamma=\bo U\bo S_0^{-1/2}$, the estimating equations for
$\bo U=(\bo u_1,\dots,\bo u_p)'$ are
\begin{equation*}
\bo u_i'\bo T(\bo u_j)=\bo u_j'\bo T(\bo u_i)\ \ \mbox{and}\ \ \bo u_i'\bo u_j=\delta_{ij},\ \ i,j=1,\dots,p,
\end{equation*}
where again $\bo T(\bo u)=\sum_{k=1}^K (\bo u'\bo R_k \bo u) \bo R_k\bo u$.
The equations then suggest a new fixed point algorithm with the two steps
\begin{align*}
&\text{step 1: } \ \   \bo T \leftarrow  (\bo T(\bo u_1),\dots,\bo T(\bo u_p))' \\
&\text{step 2:}\ \   \ \bo U\leftarrow ( \bo T \bo T')^{-1/2}  \bo T.
\end{align*}
Figure~\ref{SOBIfig} shows the estimated latent sources obtained from the data in Figure~\ref{data}. In
the signal processing literature, there are several other algorithms available for approximate simultaneous
diagonalization of $K$ matrices. The most popular one for SOBI is based on Jacobi rotations
\citep{Clarkson:1988}. Surprisingly, our new algorithm and the algorithm based on Jacobi
rotations seem to yield exactly the same solutions for practical data sets. Notice also that the symmetric
procedure does not fix the order of the rows of $\bo U$. To guarantee that
the deflation-based and  symmetric procedures estimate the same $\bo U$, we can reorder the rows of
$\bo U=(\bo u_1,\dots,\bo u_p)'$ so that
\[
\sum_{k=1}^K (\bo u_1' \bo R_k \bo u_1)^2 \ge \dots \ge \sum_{k=1}^K (\bo u_p' \bo R_k \bo u_p)^2.
\]

\section{Asymptotical properties of the SOBI estimators}
\label{properties}

In this section we derive the asymptotical properties of the two competing SOBI estimates under SOS
model~(\ref{BSSmodel}). We may assume without loss of generality that $\bo\mu=\bo 0$ and consider 
the limiting properties of the deflation-based and symmetric SOBI estimates based on the 
autocovariance matrices $\bo S_0,\bo S_1,\dots,\bo S_K$ with lags $0,1,\dots,K$. We then need some 
additional assumptions that are specific for these choices. First, we assume
\begin{itemize}
\item[(A3)]  the diagonal elements of $\sum_{k=1}^K \bo\Lambda_k^2$
are strictly decreasing.
\end{itemize}
Assumption (A3) guarantees the identifiability of the mixing matrix (with the specified autocovariance
matrices) and fixes the order of the component time series in our model. Our unmixing matrix estimate
is based on the sample autocovariance matrices $\hat{\bo S}_0,\hat{\bo S}_1,\dots,\hat{\bo S}_K$. We then further assume
that the estimates of the autocovariance matrices are root-$T$ consistent, that is,
\begin{itemize}
\item [(A4)] $\bo\Omega=\bo I_p$ and
$\sqrt{T}(\hat{\bo S}_k-\bo\Lambda_k)=O_p(1)$, $k=0,1,\dots,K$ as $T\to\infty$.
\end{itemize}
Whether (A4) is true or not depends on the distribution of the latent $p$-variate time series $\bo z$.

For the estimated autocovariance matrices, write then
$$
\hat{\bo{T}}(\bo{\gamma})=\sum_{k=1}^K(\bo{\gamma}'\hat{\bo S}_k\bo{\gamma})\hat{\bo S}_k\bo{\gamma}.
$$
The deflation-based and symmetric unmixing matrix estimates are obtained when the functionals are
applied to estimated autocovariance matrices and, consequently, they solve the following estimating
equations.

\begin{definition}
\label{estimates}
The unmixing matrix estimate
$\hat{\bo\Gamma}=(\bo{\hat\gamma}_1,\dots,\bo{\hat\gamma}_p)'$ based on $\hat{\bo S}_0$ and
$\hat{\bo S}_1,\dots,\hat{\bo S}_K$ solves the estimating equations
\[
\hat{\bo{T}}(\bo{\hat\gamma}_j)=\hat{\bo S}_0(\sum_{r=1}^j\bo{\hat\gamma}_r\bo{\hat\gamma}'_r)
\hat{\bo{T}}(\bo{\hat\gamma}_j),\ \ \ j=1,\dots,p-1,\ \ \mbox{(deflation-based)}
\]
or
\[
\hat{\bo\gamma}_i'\hat{\bo T}(\hat{\bo\gamma}_j)=\hat{\bo\gamma}_j'\hat{\bo T}(\hat{\bo\gamma}_i) \ \ \ \text{and}\ \ \
\hat{\bo\gamma}_i'\hat{\bo S}_0\hat{\bo\gamma}_j=\delta_{ij},\ \ i,j=1,\dots,p\ \ \mbox{(symmetric).}
\]
\end{definition}

Using the estimating equations and assumptions (A1)-(A4), one easily derives the following results.
The first part was already proven in~\citet{MiettinenNordhausenOjaTaskinen:2014}. For the proof of
the second part, see the Appendix.

\begin{theorem}\label{maintheorem}
Under the assumptions (A1)-(A4) we have
\begin{itemize}
\item[(i)]
the deflation-based $\hat{\bo\Gamma}=(\bo{\hat\gamma}_1,\dots,\bo{\hat \gamma}_p)'\rightarrow_p\bo I_p$,
and for $j=1,\dots,p$,
\begin{align*}
&\sqrt{T}\hat\gamma_{ji}=-\sqrt{T}\hat\gamma_{ij}-(\sqrt{T}\hat{\bo S}_0)_{ij}+o_p(1), &\ i<j, \\
&\sqrt{T}(\hat\gamma_{jj}-1)=-\frac{1}{2}\sqrt{T}((\hat{\bo S}_0)_{jj}-1)+o_p(1), &\ i=j,  \\
&\sqrt{T}\hat\gamma_{ji}=\frac{\sum_k\lambda_{kj}\left[(\sqrt{T}\hat{\bo S}_k)_{ji}
-\lambda_{kj}(\sqrt{T}\hat{\bo S}_0)_{ji}\right]}{\sum_k\lambda_{kj}(\lambda_{kj}-\lambda_{ki})}+o_p(1),
&\ i>j,
\end{align*}
\item[(ii)]
the symmetric $\hat{\bo\Gamma} =(\bo{\hat\gamma}_1,\dots,\bo{\hat\gamma}_p)' \rightarrow_p \bo I_p$, and for $j=1,\dots,p$,
\begin{align*}
\sqrt{T}\hat\gamma_{jj}&=-\frac{1}{2}\sqrt{T}((\hat{\bo S}_0)_{jj}-1)+o_p(1), &\ i=j \\
\sqrt{T}\hat\gamma_{ji}&=\frac{\sum_k(\lambda_{kj}-\lambda_{ki}) \left[(\sqrt{T}\hat{\bo S}_k)_{ji}
-\lambda_{kj}(\sqrt{T}\hat{\bo S}_0)_{ji}\right]   }{\sum_k
(\lambda_{kj}-\lambda_{ki})^2}+o_p(1), &\ i\ne j.
\end{align*}
\end{itemize}
\end{theorem}

First note that, for $\bo\Omega=\bo I_p$, the limiting distribution of the diagonal element
$\sqrt{T}(\hat\gamma_{jj}-1)$ only depends on the limiting distribution of $\sqrt{T}((\hat{\bo S}_0)_{jj}-1)$,
$j=1,\dots,p$. Hence, the comparison of the estimates should be made only using the off-diagonal
elements. Also,
\[
\sqrt{T}(\hat{\bo\Gamma}+\hat{\bo\Gamma}'-2\bo I_p) = - \sqrt{T}(\hat{\bo S}_0-\bo I_p)+o_p(1),
\]
and the limiting behavior of $\sqrt{T}(\hat{\bo\Gamma}+\hat{\bo\Gamma}'-2\bo I_p)$ is therefore similar for both
approaches.
If the joint limiting distribution of the (vectorized) autocovariance matrices is multivariate normal,
Slutsky's theorem implies that the same is true also for the unmixing matrix estimates.

\begin{corollary}
\label{asympn}
Under the assumptions (A1)-(A4), if the joint limiting distribution of
$$
\sqrt{T}\left[{\rm vec}(\hat{\bo S}_0,\hat{\bo S}_1,\dots,\hat{\bo S}_K)-{\rm
vec}(\bo I_p,\bo\Lambda_1,\dots,\bo\Lambda_K)\right]
$$
is a (singular) $(K+1)p^2$-variate normal distribution with mean value zero,
then the joint limiting distribution of $\sqrt{T}{\rm vec}(\hat{\bo\Gamma}-\bo\Gamma)$
is a singular $p^2$-variate normal distribution.
\end{corollary}

In Section~\ref{MAinf} we consider multivariate $MA(\infty)$ processes since their autocovariance
matrices have limiting joint multivariate normal distribution. So far, we have assumed that the true
value of $\bo\Omega$ is $\bo I_p$. Due to the affine equivariance of $\hat{\bo\Gamma}$, the limiting distribution
of $\sqrt{T} {\rm vec} (\hat{\bo\Gamma}\bo\Omega-\bo I_p)$ does not depend on $\bo\Omega$. 
If, for $\bo\Omega=\bo I_p$, $$\sqrt{T} {\rm vec} (\hat{\bo\Gamma}-\bo I_p)\to_d N_{p^2}(\bo 0,\bo\Sigma),$$
then, for any full-rank true $\bo\Omega$, $\hat{\bo\Gamma}-\bo\Gamma=(\hat{\bo\Gamma}\bo\Omega-\bo I_p)\bo\Gamma$
and
 $$\sqrt{T}{\rm vec} (\hat{\bo\Gamma}-\bo\Gamma)\to_d N_{p^2}\left(0,(\bo\Gamma'\otimes \bo I_p)\bo\Sigma(\bo\Gamma\otimes \bo I_p)\right). $$
Moreover, for any true $\bo\Omega$ and  $\hat{\bo\Omega}=\hat{\bo\Gamma}^{-1}$,
\[
\bo 0=\sqrt{T}\left(\hat{\bo\Gamma}\bo\Omega \bo\Gamma\hat{\bo\Omega}-\bo I_p \right)
=\sqrt{T}\,(\hat{\bo\Gamma}\bo\Omega-\bo I_p)+\sqrt{T}(\bo\Gamma\hat{\bo\Omega}-\bo I_p ) + o_P(1),
\]
which implies that
$$
\sqrt{T}{\rm vec} (\hat{\bo\Omega}-\bo\Omega)\to_d N_{p^2}\left(0,(\bo I_p\otimes\bo\Omega)\bo\Sigma(\bo I_p\otimes \bo\Omega')\right).
$$

\section{An example: MA$(\infty)$ processes}\label{MAinf}

\subsection{$MA(\infty)$ model}

An example of multivariate time series having a limiting multivariate normal
distribution is a MA$(\infty)$ process. From now on we assume that $\bo z_t$ are uncorrelated multivariate
MA$(\infty)$ processes, that is,
\begin{align}
\label{process}
\bo z_t=\sum_{j=-\infty}^\infty \Psi_j\bo \epsilon_{t-j},
\end{align}
where $\bo \epsilon_t$ are standardized iid $p$-vectors and  $\Psi_j$, $j=0,\pm 1,\pm 2,\dots$, are diagonal matrices
%with diagonal elements $\psi_{j1},\dots,\psi_{jp}$
satisfying
$\sum_{j=-\infty}^\infty \Psi_j^2=\bo I_p$. Hence
\begin{align}
\label{MAmodel}
\bo x_t=\bo\Omega \bo z_t =\sum_{j=-\infty}^\infty (\bo\Omega\Psi_j)\bo \epsilon_{t-j}
\end{align}
is also a multivariate MA$(\infty)$ process. Notice that every second-order stationary process is either a
linear process (MA($\infty$)) or can be transformed to a such one using Wold's decomposition. Notice also
that causal  ARMA$(p,q)$ processes are MA$(\infty)$ processes (see for example Chapter~3 in
\citep{BrockwellDavis:1991}).

For our assumptions, we need the following notation and definitions. We say that a $p\times p$ matrix $\bo J$
is a sign-change matrix if it is a diagonal matrix with diagonal elements $\pm 1$, and  $\bo P$ is a
$p\times p$ permutation matrix if it is obtained from an identity matrix by permuting its rows and/or
columns. For the iid $\bo \epsilon_t$, we then assume that
\begin{itemize}
\item [(B1)] $\bo\epsilon_t$ are iid with $E(\bo{\epsilon}_{t})=0$ and $Cov(\bo{\epsilon}_{t})=\bo I_p$
and with finite fourth order moments, and
\item[(B2)] the components of $\bo{\epsilon}_t$ are exchangeable and marginally symmetric, that is,
$\bo  J\bo P\bo{\epsilon}_t\sim\bo{\epsilon}_t$ for all sign-change matrices $\bo J$ and for all permutation
matrices $\bo P$.
\end{itemize}
Assumption (B1) implies that $E(\epsilon_{ti}^3\epsilon_{tj})=0$ and that
$E(\epsilon_{ti}^4)=\beta_{ii}$ and $E(\epsilon_{ti}^2\epsilon_{tj}^2)=\beta_{ij}$ are bounded
for all $i,j=1,\dots,p$. The above assumptions also imply that the model~(\ref{MAmodel})
satisfies assumptions (A1)-(A2).

\subsection{Limiting distributions of the SOBI estimates}
\label{limSOBI}

To obtain the limiting distributions of the (symmetrized) sample autocovariance matrices
$$
\hat{\bo S}_k=\frac 1{2(T-k)} \sum _{t=1}^{T-k} \left [ \bo x_t\bo x_{t+k}' + \bo x_{t+k}\bo x_t' \right],
$$
we define
$$
\bo F_k=\sum_{t=-\infty}^{\infty}\bo \psi_t \bo \psi_{t+k}', \ \ k=0,\pm 1,\pm 2,\dots
$$
where $\bo\psi_t=(\psi_{t1},\dots,\psi_{tp})'$ is the vector of the diagonal elements
of $\bo \Psi_t$, $t=0,\pm 1,\pm 2,\dots$. The diagonal elements of $F_k$ are the autocovariances of
the components of $\bo z$ at lag $k$. We also define the $p\times p$ matrices $\bo D_{lm}$, $l,m=0,\dots,K$,
with elements
\begin{align*}
(\bo D_{lm})_{ii}&=(\beta_{ii}-3)(\bo F_l)_{ii}(\bo F_m)_{ii} +
\sum_{k=-\infty}^\infty
\left( (\bo F_{k+l})_{ii} (\bo F_{k+m})_{ii}+ (\bo F_{k+l})_{ii} (\bo F_{k-m})_{ii}  \right), \\
(\bo D_{lm})_{ij}&=\frac{1}{2}\sum_{k=-\infty}^\infty((\bo F_{k+l-m})_{ii} (\bo F_k)_{jj}+
(\bo F_k)_{ii} (\bo F_{k+l-m})_{jj} ) \\
& \quad+(\beta_{ij}-1)(\bo F_l+\bo F_l')_{ij} (\bo F_m+\bo F_m')_{ij},\ \ i\neq j.
\end{align*}
The $ij$th element of $\bo D_{lm}$ is the limiting covariance of $(\hat{\bo S}_l)_{ij}$ and $(\hat{\bo S}_m)_{ij}$.
The following lemma is proved in \citet{MiettinenNordhausenOjaTaskinen:2012}.

\begin{lemma}
\label{MAasymp}
Assume that $(\bo x_1,\dots,\bo x_T)$ is a multivariate $MA(\infty)$ process defined
in~(\ref{process}) that satisfies (B1) and (B2). Then the joint limiting distribution of
$$
\sqrt{T}({\rm vec}(\hat{\bo S}_0,\hat{\bo S}_1,\dots,\hat{\bo S}_K)-{\rm
vec}(\bo I_p,\bo\Lambda_1,\dots,\bo\Lambda_K))
$$
is a singular $(K+1)p^2$-variate normal distribution with mean value zero and covariance matrix
\begin{align*}
%\label{Vmat}
\bo V=\bmat \bo V_{00} & \hdots & \bo V_{0K} \\
\vdots & \ddots & \vdots \\
\bo V_{K0} & \hdots & \bo V_{KK} \emat,
\end{align*}
with submatrices of the form
$$
\bo V_{lm}={\rm diag}({\rm vec}(\bo D_{lm}))(\bo K_{p,p}-\bo D_{p,p}+\bo I_{p^2})
$$
where
$$\bo K_{p,p}=\sum_i\sum_j(\bo e_i\bo e_j^T)\otimes(\bo e_j\bo e_i^T)
\ \ \mbox{and}\ \
\bo D_{p,p}=\sum_i(\bo e_i\bo e_i^T)\otimes(\bo e_i\bo e_i^T).$$
\end{lemma}

\bigskip

\begin{remark}
If we assume (B1) but replace (B2) by
\begin{itemize}
\item [(B2*)] the components of $\bo \epsilon_t$ are mutually independent,
\end{itemize}
then, in this independent component model case, the joint limiting
distribution of $(\hat{\bo S}_0, \hat{\bo S}_1,\dots,\hat{\bo S}_K)$ is again as given in
Lemma~\ref{MAasymp} but with $\beta_{ij} = 1$ for $i \neq j$. If we further assume that
innovations $\bo \epsilon_t$ are iid from $N_p(\bo 0, \bo I_p)$, then $\beta_{ii} = 3$
and $\beta_{ij} = 1$ for all $i \neq j$, and the variances and covariances in Lemma~\ref{MAasymp} become
even more simplified.
\end{remark}

The first part of the next theorem was presented in \cite{MiettinenNordhausenOjaTaskinen:2014},
the second part is new.

\begin{theorem}
\label{MAasympd}
Assume that $(\bo x_1,\dots,\bo x_T)$ is an observed time series from the
$MA(\infty)$ process (\ref{MAmodel})  that satisfies (B1), (B2) and (A3). Assume
(wlog) that $\bo\Omega=\bo I_p$.
If $\hat{\bo\Gamma}=(\hat{\bo\gamma}_1,\dots,\hat{\bo\gamma}_p)'$ is the SOBI estimate, then the
limiting distribution of $\sqrt{T}(\hat{\bo{\gamma}}_j-\bo e_j)$ is a $p$-variate normal distribution
with mean zero and covariance matrix
\begin{itemize}
\item [(i)] (deflation-based case)
$$
ASV(\hat{\bo{\gamma}}_j)=\sum_{r=1}^{j-1}ASV(\hat{{\gamma}}_{jr})
\bo e_r\bo e'_r+ASV(\hat{{\gamma}}_{jj})\bo e_j\bo e'_j+
\sum_{t=j+1}^pASV({\hat{\gamma}}_{jt})\bo e_t\bo e'_t,
$$
where
\begin{align*}
ASV(\hat{{\gamma}}_{jj})&=\frac 14 (\bo D_{00})_{jj}, \\
ASV(\hat{{\gamma}}_{ji})&=\frac{\mathop{\sum}_{l,m}\lambda_{li}\lambda_{mi}(\bo D_{lm})_{ji}
-2\mu_{ij}\sum_k \lambda_{ki}(\bo D_{k0})_{ji} +\mu_{ij}^2(\bo D_{00})_{ji}}{(\mu_{ij}-\mu_{ii})^2}, \\
& {\rm for}\ i<j \\
ASV(\hat{{\gamma}}_{ji})&=\frac{\mathop{\sum}_{l,m}\lambda_{lj}\lambda_{mj}(\bo D_{lm})_{ji}
-2\mu_{jj}\sum_k \lambda_{kj}(\bo D_{k0})_{ji} +\mu_{jj}^2(\bo D_{00})_{ji}}{(\mu_{jj}-\mu_{ji})^2}, \\
& {\rm for}\ i>j
\end{align*}
with $\mu_{ij}=\sum_k \lambda_{ki}\lambda_{kj}$, or
\item[(ii)] (symmetric case)
$$
ASV(\hat{\bo{\gamma}}_j)=\sum_{r=1}^{p}ASV(\hat{{\gamma}}_{jr})\bo e_r\bo e'_r
$$
where, for $i\neq j$,
\begin{align*}
ASV(\hat{{\gamma}}_{jj})=& \frac 14 (\bo D_{00})_{jj}, \\
ASV(\hat{{\gamma}}_{ji})=&\frac{\mathop{\sum}_{l,m}(\lambda_{lj}-\lambda_{li})(\lambda_{mj}-\lambda_{mi})(\bo D_{lm})_{ji}}{(\mathop{\sum}_k (\lambda_{kj}-\lambda_{ki})^2)^2} \\
&+\frac{-2\nu_{ji}\mathop{\sum}_k (\lambda_{kj}-\lambda_{ki})(\bo D_{k0})_{ji} +\nu_{ji}^2(\bo D_{00})_{ji}}{(\mathop{\sum}_k (\lambda_{kj}-\lambda_{ki})^2)^2},
\end{align*}
with $\nu_{ji}=\sum_k(\lambda_{kj}^2-\lambda_{kj}\lambda_{ki})$.
\end{itemize}
\end{theorem}

\section{Efficiency comparisons}

\subsection{Performance indices}
\label{mdsec}

In this section we compare asymptotic and finite-sample efficiencies of the two SOBI estimates.
The performance of the estimates in simulation studies can be measured using for example the minimum
distance index (MDI)~\citep{IlmonenNordhausenOjaOllila:2010b}
\begin{align*}
%\label{md}
\hat{D}= D(\hat{\bo\Gamma}\bo\Omega)=\frac{1}{\sqrt{p-1}}\inf_{\bo C\in
\mathcal{C}}\|\bo C\hat{\bo\Gamma}\bo\Omega-\bo I_p\|
\end{align*}
where $\|\cdot\|$ is the matrix (Frobenius) norm and
\[\mathcal C=\{\bo C\ :\ \mbox{each row and column of $\bo C$ has exactly one non-zero
element.}  \}\]
The minimum distance index is invariant with respect to the change of the mixing matrix, and it is
scaled so that $0\leq\hat{D}\leq 1$. It is also surprisingly easy to compute. The smaller the
MDI-value, the better is the performance. 

From the asymptotic point of view, the most attractive property of the minimum distance index is that for
an estimate $\hat{\bo\Gamma}$ with $\sqrt{T}\,{\rm vec}(\hat{\bo\Gamma}\bo\Omega-\bo I_p)\to N_{p^2}(\bo 0,\bo\Sigma)$,
the limiting distribution of $T(p-1)\hat D^2$ is that of a weighted sum of independent chi squared
variables with the expected value
\begin{align*}
%\label{md2}
{\rm tr}\left( (\bo I_{p^2}-\bo D_{p,p}) \bo\Sigma (\bo I_{p^2}-\bo D_{p,p}) \right).
\end{align*}
Notice that tr$((\bo I_{p^2}-\bo D_{p,p}) \bo\Sigma (\bo I_{p^2}-\bo D_{p,p}))$ equals the sum of the limiting
variances of the off-diagonal elements of $\sqrt{T}\,{\rm vec}(\hat{\bo\Gamma}-\bo I_p)$ and therefore
provides a global measure of the variation of the estimate $\hat{\bo\Gamma}$.

If for example  PCA, FOBI, AMUSE, deflation-based SOBI and symmetric SOBI are used to find the latent
times series based on $\bo x$ given in Figure~\ref{data}, the minimum distance index gets the values
\[  0.914,\ \ 0.628,\ \  0.246,\ \ 0.062  \ \ \mbox{and}\ \ 0.058,   \]
respectively. As PCA and FOBI are solely based on the 3-variable marginal distribution of the observations, 
they ignore time order and temporal dependence present in data. Of course, there is no reason why 
PCA should perform well here. Similarly, FOBI can be used for independent time series only if the latent 
series have distinct kurtosis values. The failure of these two methods in this example is clearly demonstrated  
by their high MDI values. AMUSE performs better here than PCA and FOBI, but is still much worse than symmetric
SOBI. Clearly the first lag is not a good choice for the separation in this example.

Finally recall that, in the signal processing literature, several other indices have been proposed for 
the finite sample comparisons of the performance of the unmixing matrix estimates
(for an overview see for example~\citet{NordhausenOllilaOja:2011}).
One of the most popular performance indices, the Amari index \citep{amari_etal:1996}, is defined as
$$
\frac 1p \left[ \sum_{i=1}^p \frac {\sum_{j=1}^p |\bo{\hat G}_{ij}|}{\max_j |\bo{\hat G}_{ij}|} +
\sum_{j=1}^p \frac {\sum_{i=1}^p|\bo{\hat G}_{ij}|}{\max_i |\bo{\hat G}_{ij}|}\right]-2,
$$
where $\hat{\bo G}=\hat{\bo\Gamma}\bo\Omega$.
The index is invariant under permutations and sign changes of the rows and columns of $\bo{\hat G}$.
However, heterogeneous rescaling of the rows (or columns) on $\bo{\hat G}$ changes its value. Therefore,
for the comparisons, the rows of $\hat{\bo\Gamma}$ should be rescaled in a similar way. We prefer MDI,
since the Amari index is based on the $L_1$ norm and cannot be easily related to the limiting
distribution of the unmixing matrix estimate.

\subsection{Four models for the comparison}

The following four models were chosen for the comparison of the deflation-based and symmetric SOBI
estimates.
The components of the source vectors are
\begin{itemize}
\item[(a)] three MA(10)-series with coefficient vectors \\
$(0.8,3.8,1.2,1.4,1.1,0.5,0.7,0.3,0.5,1.8)$, \\
$(-0.6,1.3,-0.1,1.3,1.6,0.4,0.5,-0.4,0.1,2.8)$ and \\
$(-0.4,-1.5,0,-1.1,-1.9,0,-0.7,-0.4,-0.2,0.4)$, \\
respectively, and normal innovations,
\item[(b)] three AR-series with coefficient vectors $(0.6)$, $(0,0.6)$ and
$(0,0,0.6)$, respectively, and normal innovations,
\item[(c)] three ARMA-series with AR-coefficient vectors $(0.3,0.3,-0.4)$, \\
$(0.2,0.1,-0.4)$, $(0.2,0.2,0.4)$, and  MA-coefficient vectors \\
$(-0.6,0.3,1.1,1.0,-1.1,-0.3)$, $(1.2,2.8,-1.0,-1.0,0.1,0.1)$, \\
$(-1.4,-1.9,-0.5,-0.3,-0.4,0.4)$, respectively, and normal innovations,
\item[(d)] three AR(1)-series with coefficients $0.6$, $0.4$ and $0.2$, respectively, and normal
innovations.
\end{itemize}
Each component is scaled to unit variance. Due to the affine equivariance of the estimates, it is not a
restriction to use $\bo\Omega=\bo I_3$ in the comparisons.

\subsection{Asymptotic efficiency}

The asymptotic efficiency of the estimates can be compared using the sum of the limiting variances
of the off-diagonal elements of $\sqrt{T}\,{\rm vec}(\hat{\bo\Gamma}-\bo I_p)$. In Table~\ref{table1}, these
values are listed for the symmetric and deflation-based SOBI estimates, when both methods are using
lags $k=1,\dots,10$ in all four models.

\begin{table}[htb]
\caption{\label{table1}The sums of the estimated variances of selected rows of
$\hat{\bo\Gamma}$ for the symmetric SOBI estimates utilizing four candidate sets of lags.}
\centering
   \begin{tabular}{lcccc}
    \hline
     &  \multicolumn{4}{c}{model} \\ \cline{2-5}
     & (a) & (b) & (c) & (d) \\
    \hline
     deflation-based & 46.5 & 31.8 & 11.0 & 61.6  \\
     symmetric & 24.1 & 10.6 & 9.4 & 75.1 \\
    \hline
    \end{tabular}
\end{table}

First notice that model (d) is the only one, where the deflation-based  estimate outperforms the
symmetric estimate. Also in model (c) the deflation-based method has quite competitive variances whereas
in models (a) and (b) the symmetric estimate is much more accurate than the deflation-based
estimate.

\subsection{Finite sample behavior}

\begin{figure}[htb]
\centering
\includegraphics[width=12truecm,height=12truecm]{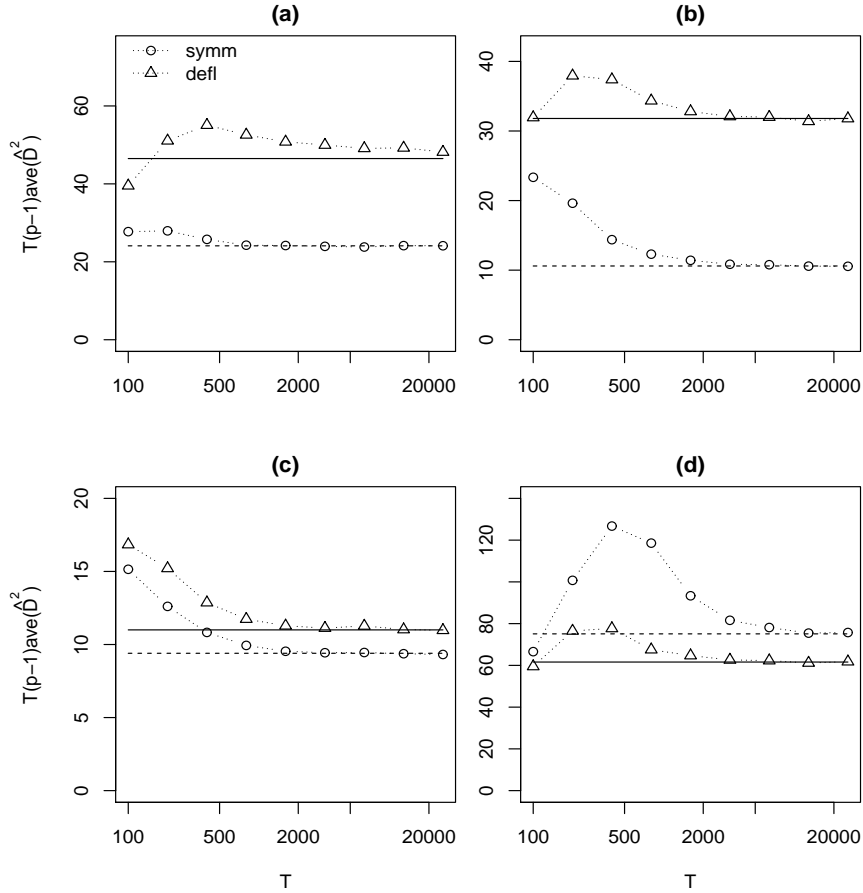}
\caption{\label{fig1}The averages of $T(p-1)\hat{D^2}$ for the symmetric and deflation-based SOBI
estimates from 10 000 repetitions of observed time series with length $T$ from models (a)-(d).
The horizontal lines give the expected values of the limiting distributions of $T(p-1)\hat{D}^2$. }
\end{figure}

For finite sample efficiencies, we use simulated time series from models (a)-(d). For each model and for
each value of $T$, we have 10 000 repetitions of the simulated time series. The averages of
$T(p-1)\hat{D}^2$ were then computed for the symmetric SOBI and deflation-based SOBI estimates.
Again, we use lags $k=1,\dots,10$ for both estimates. The results are plotted in Figure~\ref{fig1}.

As explained  in Section~\ref{mdsec}, the limiting distribution of  $T(p-1)\hat{D}^2$ has expected
values as given in Table~\ref{table1}. We then expect that the averages of $T(p-1)\hat{D}^2$ converge
to these expected values in all four models. As seen in Figure~\ref{fig1}, this seems to happen
but the convergence is quite slow in all cases. Notice also that the finite sample efficiencies are 
higher than their approximations given by the limiting results. 

Functions to compute deflation-based and 
symmetric SOBI estimates and their theoretic and finite sample
efficiencies are provided in the R packages JADE \citep{JADEr:2013} and BSSasymp
\citep{MiettinenNordhausenOjaTaskinen:2013}.

\section{EEG example}
\label{EEG}

\begin{figure}[htb]
\centering
\includegraphics[width=12truecm,height=5truecm]{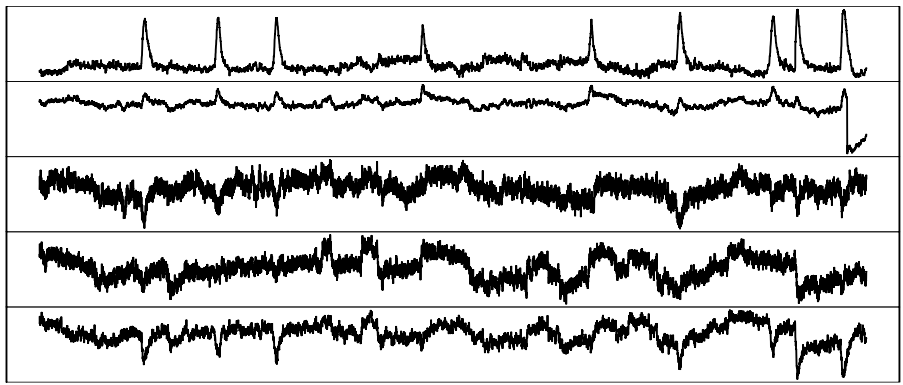}
\caption{\label{EEGsignal} Five randomly selected EEG components of length $T=10000$, from 20001 to 30000.}
\end{figure}

Let us now illustrate how the theoretical results derived in Section~\ref{limSOBI} can be used to select 
a good set of autocovariance matrices to be used in SOBI method for a problem at hand. We consider EEG 
(electroencephalography) data recorded at the Department of Psychology, University of Jyv\"askyl\"a, from an 
adult using 129 electrodes placed on the scalp and face. The EEG data was online bandpass filtered at 0.1 - 100 Hz 
and sampled at 500 Hz. The length of the data set equals $T=100000$. An extract of five randomly selected EEG components is plotted in Figure~\ref{EEGsignal}.

The goal in EEG data analysis is to measure the brain's electrical activity. As the measurements are made along the 
scalp and face, a mixture of unknown source signals is observed. Moreover, EEG recordings are often contaminated
by non-cerebral artifacts, such as eye movements and muscular activity. Our aim is to use SOBI methods to 
separate few clear artifact components from EEG data. When choosing the set of autocovariance matrices to be 
diagonalized in SOBI, we follow the recommendations of~\cite{TangLiuSutherland:2005} and compare the results 
given by the following four sets of lags.
\begin{itemize}
\item[(1)] $1,2,\dots,10,12,\dots,20,25,\dots,100,120,\dots,300$.
\item[(2)] $1,2,\dots,10,12,\dots,20$.
\item[(3)] $1,2,\dots,10,12,\dots,20,25,\dots,100$.
\item[(4)] $25,30,\dots,100,120,\dots,300$.
\end{itemize}

We applied both, deflation-based and symmetric, SOBI methods to the EEG data using the four different 
autocovariance sets. However, as the deflation-based method
gave significantly larger variance estimates for the unmixing matrix estimate, we only report the results
based on symmetric SOBI. We chose three recognizable components, i.e. eye blink, horizontal
eye movement and muscle activity, and estimated the sum of variances of the corresponding rows of the unmixing 
matrix under the assumption that the source signals are generated by $MA(\infty)$ processes with normal innovations.
The resulting variance estimates are reported in Table~\ref{table2}.

\begin{table}[htb]
\caption{\label{table2}The sums of the estimated variances of selected rows of
$\hat{\bo\Gamma}$ for the symmetric SOBI estimates utilizing four candidate sets of lags.}
\centering
   \begin{tabular}{lcccc}
    \hline
     &  \multicolumn{4}{c}{set of lags} \\ \cline{2-5}
     & (1) & (2) & (3) & (4) \\
    \hline
     eye blink & 0.052 & 0.021 & 0.060 & 0.091  \\
     horizontal eye movement & 0.234 & 0.192 & 0.337 & 0.517 \\
     muscle activity & 0.058 & 0.051 &  0.057 & 0.123    \\
    \hline
    \end{tabular}
\end{table}

%\begin{table}
%\caption{\label{table3}The sums of the estimated variances of selected rows of
%$\hat{\bo\Gamma}$ for the deflation-based SOBI estimates.}
%\centering
%\makebox{%
%   \begin{tabular}{*{3}{c}}
%    \hline
%    & blink & muscle  \\
%    \hline
%     (1) & 0.579 & 0.107 \\
%     (2) & 5.484 & 0.071  \\
%     (3) & 1.002 & 0.101  \\
%     (4) & 2.884 & 1.751  \\
%    \hline
%    \end{tabular}}
%\end{table}

\begin{figure}[htbp]
\centering
\includegraphics[width=12truecm,height=4truecm]{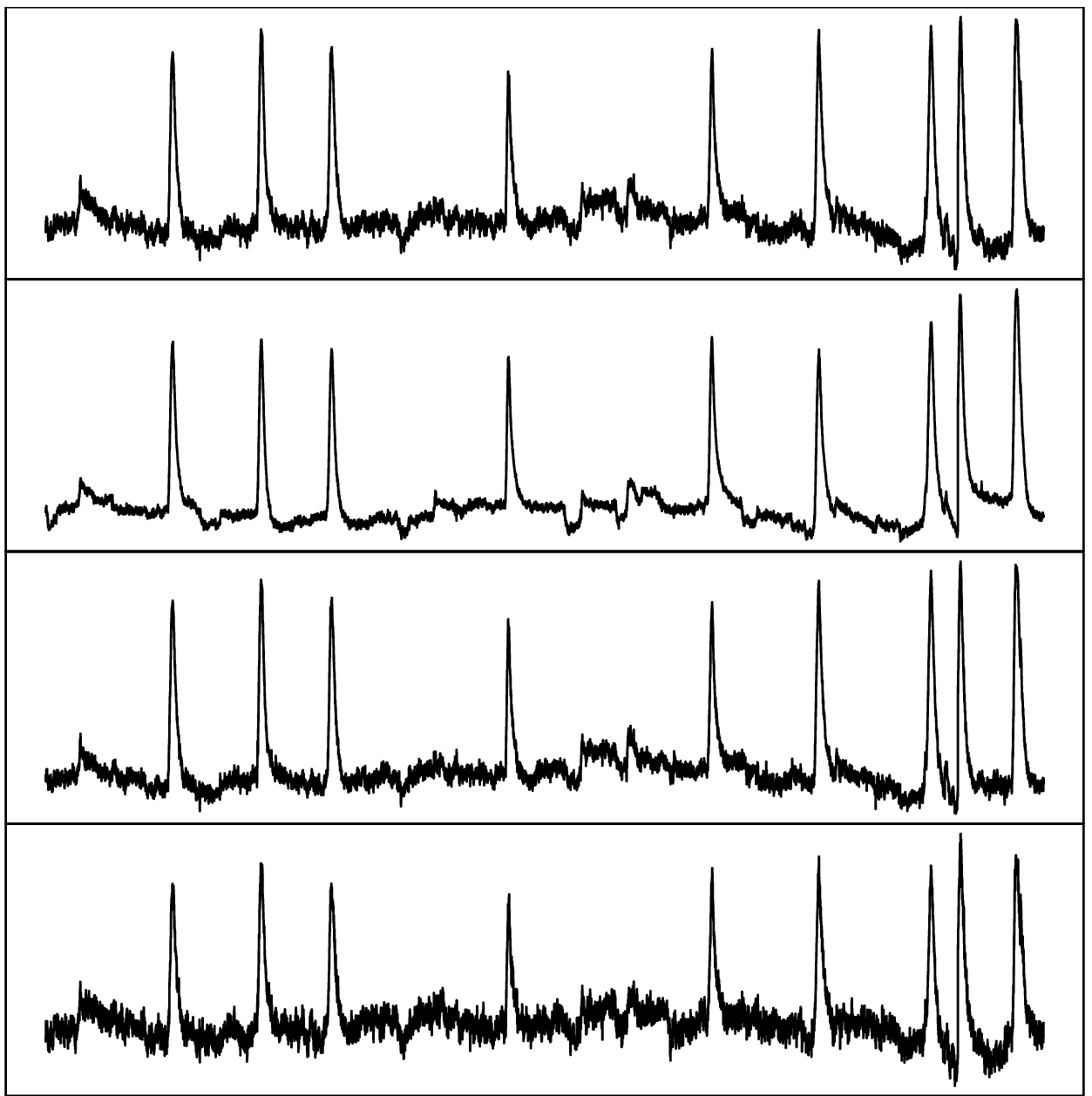}
\vskip0.5truecm
\includegraphics[width=12truecm,height=4truecm]{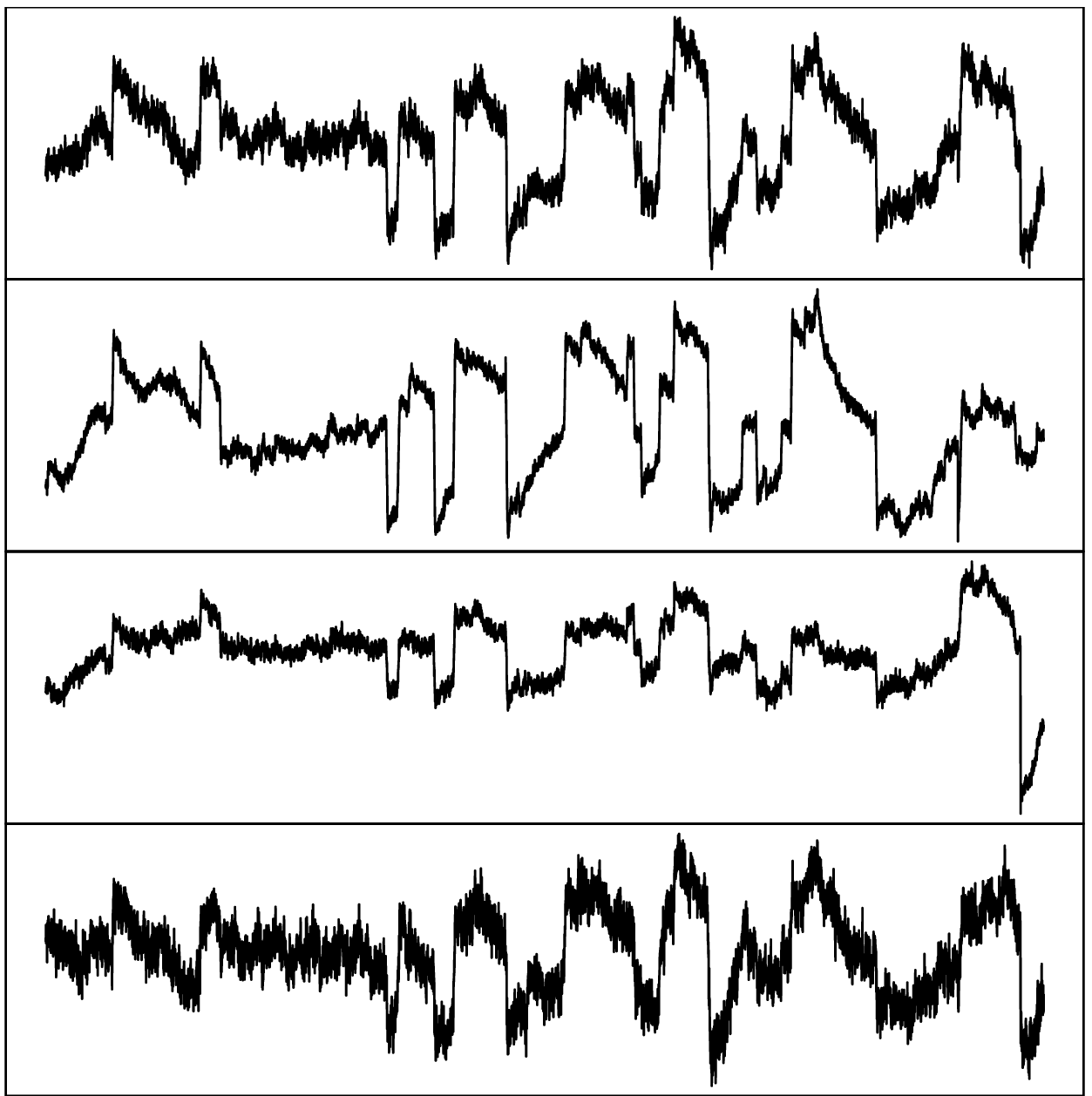}
\vskip0.5truecm
\includegraphics[width=12truecm,height=4truecm]{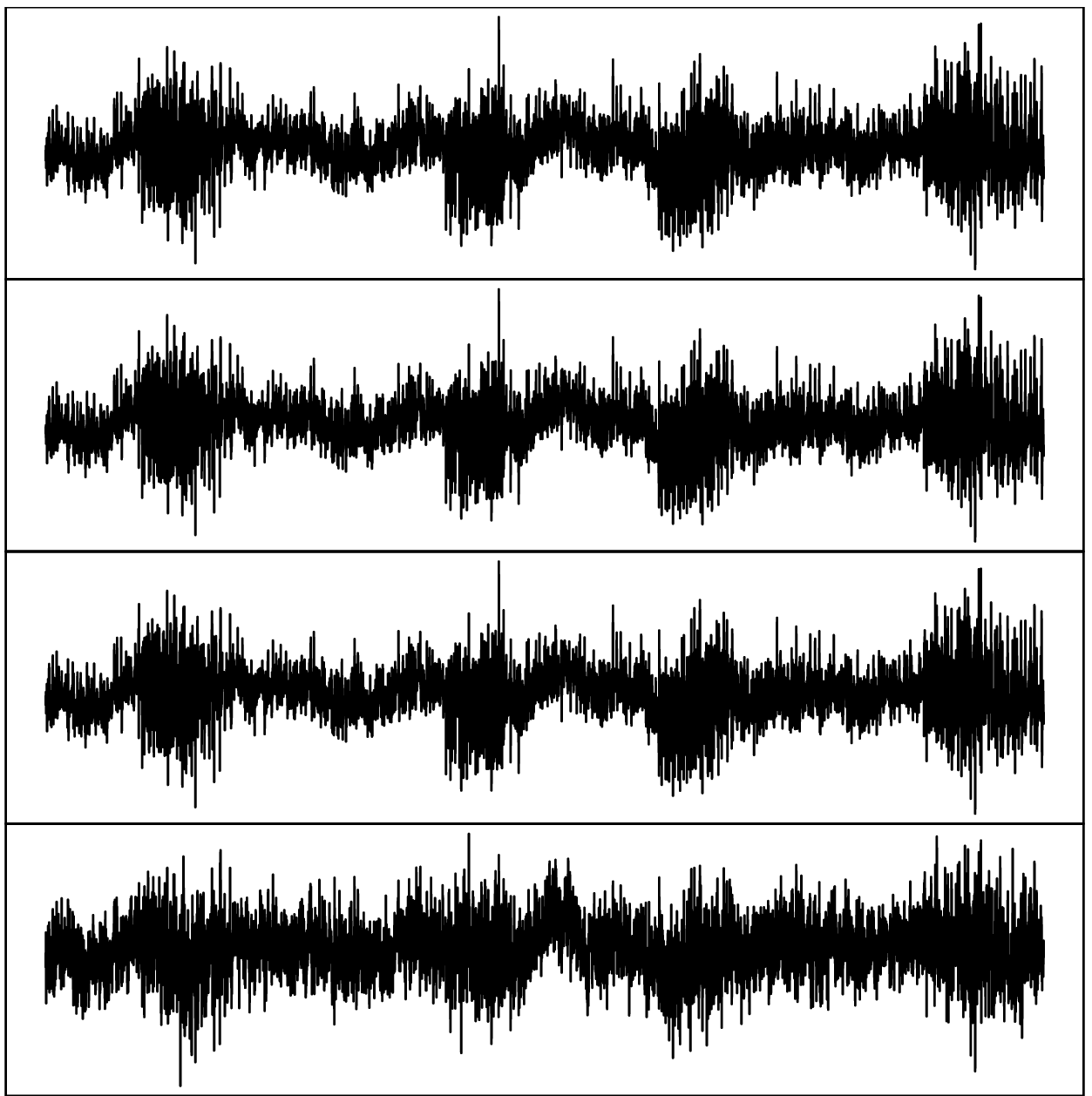}
\caption{\label{EEGest} Estimated eye blink (top figure), horizontal eye movement (middle figure) and 
muscle activity (bottom figure) components of length $T=10000$, from 20001 to 30000. The components 
are estimated using symmetric SOBI method with lags given in sets (1)-(4), respectively.}
\end{figure}

%In the study reported by~\cite{TangLiuSutherland:2005}, the SOBI method that was based on autocovariance
%matrices with lags in set (1) was shown to have the best overall performance among its competitors. In our example,
As the results of Table~\ref{table2} indicate, the symmetric SOBI method that used autocovariance 
matrices with lags in set (2) gave best separation results for the three components of interest. 
These findings are confirmed by the time plots of separated components in Figure~\ref{EEGest}. The 
good performance of symmetric SOBI based on lags in set (2) is especially visible when looking at 
the corresponding eye blink and horizontal eye movement components.

%\begin{figure}[htbp]
%\centering
%\makebox{\includegraphics[width=6truecm,height=6truecm]{EEG_blink_def.eps}
%\includegraphics[width=6truecm,height=6truecm]{EEG_muscle_def.eps}}
%\caption{\label{EEGest3} Blink and muscle activity components estimated with deflation-based (1)-(4). Length is 100%00.}
%\end{figure}

\section{Concluding remarks}
\label{conclusion}

Symmetric SOBI is a popular blind source separation method but a careful analysis of its statistical
properties has been missing so far. The theoretical results for deflation-based SOBI were presented only
recently in \cite{MiettinenNordhausenOjaTaskinen:2014}.
There is a lot of empirical evidence that symmetric BSS methods perform better than their
deflation-based competitors. In this paper we used the Lagrange multiplier technique to derive
estimating equations for symmetric SOBI that allowed a thorough theoretic analysis of the
properties of the estimate. In most cases we studied, the limiting efficiency of the symmetric
SOBI estimate was better than that of the deflation-based estimate. The estimating equations also
suggested a new algorithm for symmetric SOBI. Such an algorithm gave, in all our simulations, exactly 
the same results as the most popular algorithm based on Jacobi rotations. In a separate paper, these and
other algorithms with associated estimates are compared in various settings with different values
of $p$ and $K$.

The problem corresponding to the selection of lags is still open; only few ad-hoc
guidelines are available in the literature, see for example \cite{TangLiuSutherland:2005}.
We followed these guidelines in our EEG data example. Notice however that the results presented in this
paper can be used to build a two-stage estimation procedure where, at stage 2, the final SOBI estimate
is selected among all SOBI estimates using their estimated efficiencies in a model determined
by a preliminary SOBI estimate applied at stage 1. The results derived here can also be applied 
to different inference procedures, including hypothesis testing and model selection.

\section*{Acknowledgements}
This research was supported by the Academy of Finland (grants 256291 and 268703). We thank
Dr. Jarmo H\"am\"al\"ainen for providing us with the EEG data.

\section*{Appendix}
\label{AppA}

\noindent{\bf Proof of Theorem~\ref{maintheorem}}\ \
\\

\noindent(i) The proof for the consistency and limiting behavior of the deflation-based SOBI estimate can be found in
\cite{MiettinenNordhausenOjaTaskinen:2014}.
\bigskip

\noindent(ii) We first prove the consistency of the estimate.
Let $\mathcal{\bo U}$ be the compact set of all $p\times p$ orthogonal matrices.
For $\bo U=(\bo u_1,\dots,\bo u_p)'\in\mathcal{\bo U}$, write
\[
\Delta(\bo U)=\sum_{i=1}^p\sum_{k=1}^K (\bo u_i'\bo R_k\bo u_i)^2 \ \ \mbox{and}\ \
\hat{\Delta}(\bo U)=\sum_{i=1}^p\sum_{k=1}^K (\bo u_i'\bo{\hat R}_k\bo u_i)^2
\]
As
$(\bo u'\bo{\hat R}_k \bo u)^2-(\bo u'\hat R\bo u)^2=(\bo u'(\bo{\hat R}_k-\bo R_k)\bo u)(\bo u'(\bo{\hat R}_k+\bo R_k)\bo u)
$
and
$\bo{\hat R}_k\to_P \bo R_k$, for $k=1,\dots,K$,
\[
M=\sup_{\bo U\in\mathcal{\bo U}} |\Delta(\bo U)-\hat{\Delta}(\bo U)|\le 2\sum_{k=1}^K ||\bo{\hat R}_k-\bo R_k||\to_P 0.
\]
Under our assumptions, $\bo U=\bo I_p$ is the unique maximizer  of $\Delta(\bo U)$ in the subspace
\[
\mathcal{\bo U}_0=\left\{ \bo U\in \mathcal{\bo U}\ :\  \mbox{$\bo u_i\bo R_i\bo u_i$ descending and $\bo u_i'\bo 1_p\ge 0$, $i=1,\dots,p$}\right\}.
\]
Notice that, in this subspace $\mathcal{\bo U}$, the order and signs of the rows of $\bo U$ are fixed.
For all $\epsilon>0$, write next
\[
\mathcal{\bo U}_\epsilon=\left\{ \bo U\in \mathcal{\bo U}_0\ :\ ||\bo U-\bo I_p||\ge \epsilon \right\}
\]
and
\[
\delta_\epsilon=\Delta(\bo I_p)-\sup_{\bo U\in \mathcal{\bo U}_\epsilon} \Delta(\bo U).
\]
Clearly, $\delta_\epsilon>0$ and $\delta_\epsilon\to 0$ as $\epsilon\to 0$.
Let $\hat{\bo U}$ be the unique maximizer of $\hat{\Delta}(\bo U)$ in $\mathcal{\bo U}_0$.
Then
\begin{eqnarray*}
%P\left(||\hat{\bo U}-\bo I_p||<\epsilon\right) & \ge & P\left(\hat{\Delta}(\bo I_p)>
%\sup_{\bo U\in \mathcal{\bo U}_\epsilon}\hat{\Delta}(\bo U)\right) \\
%  \ & \ge & P\left(M\le \frac {\delta_\epsilon} 3\right)\to 1, \ \ \forall \epsilon
P\left(||\hat{\bo U}-\bo I_p||<\epsilon\right)  \ge  P\left(\hat{\Delta}(\bo I_p)>
\sup_{\bo U\in \mathcal{\bo U}_\epsilon}\hat{\Delta}(\bo U)\right)
\ge P\left(M\le \delta_\epsilon / 3\right)\to 1, \ \ \forall \epsilon
\end{eqnarray*}
and the convergence $\hat{\bo U}\to_P \bo I_p$ follows. Thus
$\hat{\bo\Gamma}=\hat{\bo U}\hat{\bo S}_0^{-1/2}\to_P \bo I_p$ also holds true.

To prove the second part of the result (ii), notice that the estimating equations give
\begin{equation}\label{eq1}
0=\sqrt{T}\bo{\hat\gamma}_i'\hat{\bo S}_0\bo{\hat\gamma}_j=\sqrt{T}\hat\gamma_{ij}+\sqrt{T}\hat \gamma_{ji}+\sqrt{T}(\hat{\bo S}_0)_{ij}+o_p(1),
\end{equation}
for $i\ne j$, and
\begin{equation}\label{eq2}
0=\sqrt{T}(\bo{\hat\gamma}_i'\hat{\bo S}_0\bo{\hat\gamma}_i-1)=2\sqrt{T}(\hat\gamma_{ii}-1)+\sqrt{T}((\hat{\bo S}_0)_{ii}-1)+o_p(1).
\end{equation}
Next note that
\begin{equation}\label{eq3}
\hat{\bo\gamma}_i'\bo{\hat T}(\hat{\bo\gamma}_j)-\bo e_i'\bo T(\bo e_j)=\hat{\bo\gamma}_j'\bo{\hat T}(\hat{\bo\gamma}_i)-\bo e_j'\bo T(\bo e_i)\end{equation}
and
\begin{equation}\label{eq4}
\begin{split}
\sqrt{n}\left(\hat{\bo\gamma}_i'\bo{\hat T}(\hat{\bo\gamma}_j)-\bo e_i'\bo T(\bo e_j)\right)& =\sum_{k=1}^K
\lambda_{kj}^2\sqrt{n}\hat \gamma_{ij}+ \sum_{k=1}^K \lambda_{kj}\sqrt{n}(\hat{\bo S}_k)_{ij} \\
 & \quad + \sum_{k=1}^K \lambda_{ki}\lambda_{kj}\sqrt{n}\hat\gamma_{ji}+o_p(1).
\end{split}
\end{equation}
for all $i\ne j$. The result then follows from equations (\ref{eq1})-(\ref{eq4}).
\end{document}